\documentclass[11pt,a4paper]{amsart}
\usepackage[utf8]{inputenc}
\usepackage{amsmath}
\usepackage{amsfonts}
\usepackage{amssymb}
\usepackage{graphicx}
\usepackage{amsmath,amssymb, amsthm}
\usepackage{pgf}
\usepackage{tikz-cd}
\usepackage{xy}
\usepackage{hyperref}
\usetikzlibrary{matrix,arrows,decorations.pathmorphing}
\usepackage[english]{babel}
\usepackage{mathtools}
\usepackage{caption}
\usepackage{geometry}
\usepackage{float}

\theoremstyle{definition}

\theoremstyle{definition}
\newtheorem{rmk}{Remark}
\theoremstyle{definition}

\theoremstyle{definition}\usepackage{amsmath}

\theoremstyle{plain}
\newtheorem{thm}{Theorem}
\theoremstyle{plain}

\theoremstyle{plain}

\theoremstyle{plain}

\newcommand{\Z}{\mathbb{Z}}
\newcommand{\Q}{\mathbb{Q}}

\newcommand{\Ok}{\mathcal{O}_K}

\newcommand{\Tr}{\operatorname{Tr}}
\newcommand{\Nm}{\operatorname{N}}

\newcommand{\lp}{\left(}
\newcommand{\rp}{\right)}

\newcommand{\pp}{\mathfrak{p}}

\begin{document}

\address{Dipartimento di Matematica\\
         Università di Milano\\
         via Saldini 50\\
         20133 Milano\\
         Italy}
\title[On small discriminants of number fields of degree 8 and 9]{On small discriminants of number fields of degree 8 and 9}
\author[F.~Battistoni]{Francesco Battistoni}

\email{francesco.battistoni@unimi.it}



\maketitle

\begin{abstract}
    We classify all the number fields with signature (4,2), (6,1), (1,4) and (3,3) having discriminant lower than a specific upper bound. This completes the search for minimum discriminants for fields of degree 8 and continues it in the degree 9 case. We recall the theoretical tools and the algorithmic steps upon which our procedure is based, then we focus on the novelties due to a new implementation of this process on the computer algebra system PARI/GP; finally, we make some remarks about the final results, among which the existence of a number field with signature $(3,3)$ and small discriminant which was not previously known.
\end{abstract}

\section{Introduction}
Consider the family of number fields $K$ with fixed degree $n$ and fixed signature $(r_1,r_2)$. Classical results by Minkowski and Hermite, obtained at the end of the XIX-th century, imply the following properties for the discriminant $d_K$ of the fields:
\begin{itemize}
    \item There exists an explicit lower bound $|d_K| > C(n,r_2)$, where $C(n,r_2)>1$ depends only on the degree and the signature. This implies that in the family of fields with degree $n$ and signature $(r_1,r_2)$ there exists a field $F$ such that $|d_K|$ attains the minimum value when $K=F$.
    \item For every $C>0$, there exist only finitely many number fields with fixed degree and signature and such that $|d_K|\leq C$.
\end{itemize}
The study of number fields with respect to their discriminants is then characterized by a double purpose: to find the minimum values for the discriminants of fields with fixed signature, and to completely classify all the fields in this family up to a chosen discriminant bound (a goal which encompasses the first one). Complete tables of number fields with bounded discriminant are useful tools in Number Theory, because they provide explicit examples over which one can get some heuristic or prove results which are known to be asymptotically true in the discriminant (see, as an example, the work by Astudillo, Diaz y Diaz and Friedman \cite{regulators} on minimum regulators which explicitly requires this kind of lists).

Giving a complete classification of fields with fixed degree $n$ and signature $(r_1,r_2)$ is easy for $n=2$, because any quadratic field has the form $\Q(\sqrt{d})$ with $d\in\Z$ squarefree, and this structure returns $d_K$ equal to either $4d$ or $d$, depending on the residue class of $d$ modulo 4.

For $n=3$ the research is still not difficult, thanks to Davenport-Heilbronn's correspondance between isomorphism classes of cubic number fields and equivalence classes of primitive binary integral cubic forms \cite{davenportHeilbronn}, which was the theoretical cornerstone for Belabas' algorithm for the classification of cubic fields with bounded discriminant \cite{belabasCubic}.

Whenever one considers fields of higher degree, the classification becomes harder to get. There are two main mathematical frameworks, developed during the 70's and 80's of the XX-th century, which allowed several researchers to get results for fields with low degree:
\begin{itemize}
    \item Geometry of Numbers and its applications to the rings of integers, which provided explicit estimates on the possible maximum values for the coefficients of the defining polynomials of number fields $K$ with bounded $|d_K|$: this was investigated by Hunter and Pohst \cite{pohst1982computation} for number fields over $\Q$ and by Martinet \cite{martinet1985methodes} for generic number field extensions;
    \item Lower bounds for the discriminants derived from the explicit formulae of Dedekind Zeta functions, a procedure which was pursued by Odlyzko \cite{odlyzko1990bounds}, Poitou \cite{poitou1977petits} and Serre \cite{serre1986minorations} and which allowed Diaz y Diaz \cite{y1980tables} to obtain lower bounds of $|d_K|$ for several degrees and signatures. 
\end{itemize}
The simultaneous use of the previous tools permitted to develop algorithmic procedures which gave complete classifications of number fields up to certain discriminant bounds in the following cases:
\begin{itemize}
    \item Number fields with degree 4 \cite{QuarticFields} and 5 \cite{QuinticFields};
    \item Number fields with degree 6 \cite{bergeOlivierMartinet, olivierSexticQuadratic1, olivierQuadratic2, olivierPrimitifs, olivierSextiquesCubique, olivierSextiquesPrimitifs, olivierSexticCubic1992} and degree 7 \cite{diazSepticUneSeulePlace, diazSepticCertain, diazSeptic, pohstDegree7TotallyReal}.
    \item Totally complex \cite{y1987petits} and totally real \cite{pohst1990minimum} number fields of degree 8;
    \item Totally real number fields of degree 9 \cite{takeuchiNine}.
\end{itemize}
For what concerns further signatures in degree 8, no complete tables up to some bound were known and for several years no attempts of this kind were made. During his Ph.D. work, the author \cite{battistoniMinimum} was then able to give a complete classification of number fields with degree 8, signature $(2,3)$ and $|d_K|\leq 5726301$, showing that there exist exactly 56 such fields: this result was obtained by exploiting the aforementioned theoretical ideas in order to write an algorithmic procedure which was implemented in a program relying on the softwares MATLAB and PARI/GP \cite{pari}. 

This setting was not good enough for exploring other signatures in degree 8 and 9, but we are now able to provide a better implementation, needing just PARI/GP, which allowed us to obtain the following classification result. 
\begin{thm}\label{TeoremaClassificazioneDiscriminanti}
There exist 41 number fields $K$ with signature $(4,2)$ and with $|d_K|\leq 20829049$. The minimum value of $|d_K|$ is $15243125$.\\
There exist 8 number fields $K$ with signature $(6,1)$ and with $|d_K|\leq 79259702$. The minimum value of $|d_K|$ is $65106259$.\\
There exist 67  number fields $K$ with signature $(1,4)$ and with $|d_K|\leq 39657561$. The minimum value of $|d_K|$ is $29510281$.\\
There exist 116 number fields $K$ with signature $(3,3)$ and with $|d_K|\leq 146723910 $. The minimum value of $|d_K|$ is $109880167$.
\end{thm}
\noindent
The number fields and the complete tables are gathered in the website \cite{battistoniWebsite}, together with the PARI/GP programs used for their classification. The programs were run on the cluster system of Université de Bordeaux and on the clusters INDACO and HORIZON of Università degli Studi di Milano.
\\\\
Here is an overview of the paper.

Section 2 recalls the theoretical foundations of the algorithmic procedure, which are respectively Hunter-Pohst-Martinet's Theorem and the local corrections for lower bounds of discriminants given by prime ideals. Section 3 presents the various steps in which the algorithm is divided. Section 4 finally presents the main novelties of our new implementation and some remarks on the final results obtained.

\subsection*{Acknowledgements}
I would like to thank my Ph.D. advisor Giuseppe Molteni for every useful suggestion and for his supervision of my research in Università degli Studi di Milano. I would like to thank also Institut de Mathématiques de Bordeaux, for hosting me and allowing me to use the IMB cluster, and the people I worked with and gave me many advices: Bill Allombert, Karim Belabas, Henri Cohen, Andreas Enge, Aurel Page, Guillame Ricotta, Damien Robert.\\
Thanks also to Alessio Alessi and Francesco Fichera, who gave me permission for using INDACO and HORIZON clusters respectively, to Scherhazade Selmane, for lending me her tables with local corrections, and to Gunter Malle and Ken Yamamura, for their remarks about the preliminary version of this paper.

\section{Theoretical recalls}
In our procedure we look for irreducible monic polynomials of degree 8 and 9 with integer coefficients which define the desired number fields: the first problem consists then in giving an upper bound to the number of these polynomials, and the bound should depend on the discriminant and the signature.

Given a number field $K$, an element $\alpha\in \Ok$ and a number $k\in\Z$, the \textbf{Newton sum of order $k$ of $\alpha$} is defined as the sum
$$S_k(\alpha):= \sum_{i=1}^n \alpha_i^k$$
where the $\alpha_i$'s represent the conjugates of $\alpha$ with respects to the embeddings $\sigma_1,\ldots,\sigma_n$ of $K$.

One has $S_k(\alpha)\in\Z$, and $S_1(\alpha)=\Tr(\alpha)$; moreover, if $f(x):= x^n +a_1 x^{n-1}+\cdots+ a_{n-1}x+a_n$ is the defining polynomial of $\alpha$, then one has the recursive relations
\begin{equation}\label{recursiveRelations}
     S_k(\alpha) = -k a_k -\sum_{j=1}^{k-1} a_j S_{k-j}(\alpha)\hspace{0.5cm} \text{ for every }2\leq k\leq n
\end{equation}
which link the coefficients of $f(x)$ to the values of the Newton sums.

Consider then the absolute Newton sum
$$T_2(\alpha):= \sum_{i=1}^n |\alpha_i|^2.$$
We have an estimate for $T_2(\alpha)$, depending on $n$ and $|d_K|$, provided by Hunter-Pohst-Martinet's Theorem \cite{martinet1985methodes}.

\begin{thm}\label{HunterPohstMartinet}
    Let $K$ be a number field of degree $n$ with discriminant $d_K$. Then there exists an element $\alpha\in\Ok\setminus\Z$ which satisfies the following conditions:
\begin{align}
    &\text{A) }  0\leq \Tr(\alpha)\leq \left\lfloor\frac{n}{2} \right\rfloor; \nonumber\\
    &\text{B) }  T_2(\alpha) \leq \frac{\Tr(\alpha)^2}{n} + \gamma_{n-1}\lp\frac{|d_K|}{n}\rp^{1/(n-1)}=:U_2\nonumber
\end{align}
where $\gamma_{n-1}$ is Hermite's constant of dimension $n-1$.
\end{thm}

The element $\alpha$ is called an \textbf{HPM-element} for $K$: for such an algebraic integer, the previous theorem allows us to compute an upper bound for its trace and its second Newton sum $S_2$. These data, together with the absolute value of the norm $N:=|\Nm(\alpha)|$, are enough for giving upper bounds to every Newton sum, thanks to Pohst's result \cite{pohst1982computation}.
\begin{thm}\label{ThmPohst}
Given $K$, $U_2$ and $\alpha$ as in Theorem \ref{HunterPohstMartinet}, given $N\in\mathbb{N}$ such that $N\leq(U_2/n)^{n/2}$, then for every $k\in\Z\setminus\{0,2\}$ we have an inequality
$$|S_k(\alpha)|\leq U_k$$
where $U_k$ is a number depending on $n,r_1$ and $U_2$.
\end{thm}
Our goal is then to test the polynomials generated by a choice of the coefficients which derives from the values of the Newton sums $S_k$ (with $2\leq k\leq n$) ranging in the intervals $[-U_k,U_k]$ and satisfying the recursive relations \eqref{recursiveRelations}. In order to do so, we need to choose an upper bound for $|d_K|$.

\begin{rmk}
The condition $N\leq (U_2/n)^{n/2}$ is set in order to respect the inequality between geometric and arithmetic means: in fact,
$$N^2 = \prod_{i=1}^n |\alpha_i|^2 \leq \lp\frac{\sum_{i=1}^n |\alpha_i|^2}{n}\rp^n = \lp\frac{U_2}{n}\rp^n.$$
\end{rmk}

We recall now an inequality, proved by Poitou \cite{poitou1977petits}, which gives a lower bound for the discriminants of number fields with fixed degree $n$ and signature $(r_1,r_2)$.
\begin{thm}
Let $K$ be a number field of degree $n$, signature $(r_1,r_2)$ and discriminant $d_K$. Let $f(x)$ be the function

$$	f(x):= \lp \frac{3}{x^3}(\sin x -x\cos x)\rp^2.$$
Then, for every $y>0$, one has
\begin{equation}\label{RootDiscriminantLowerBoundNumerical}
		\frac{1}{n}\log|d_K| \geq \gamma + \log 4\pi - L_1(y) -\frac{12\pi}{5n\sqrt{y}} + \frac{4}{n}\sum_{\pp\subset\Ok}\sum_{m=1}^{\infty}\frac{\log \Nm(\pp)}{1+(\Nm(\pp))^m}f(m\sqrt{y}\log \Nm(\pp))
\end{equation}
where $\gamma$ is Euler's constant, the sum runs over the non-zero prime ideals of $\Ok$, $\Nm(\pp)$ is the absolute norm of the prime $\pp$ and

\begin{equation*}
	L_1(y):= L(y) +\frac{1}{3}L\lp\frac{y}{3^2}\rp + \frac{1}{5}L\lp\frac{y}{5^2}\rp +\cdots +\frac{r_1}{n}\left[L(y)-L\lp\frac{y}{2^2}\rp+L\lp\frac{y}{3^2}\rp\cdots\right]
\end{equation*}
where

\begin{equation*}
L(y):=-\frac{3}{20y^2} + \frac{33}{10y} + 2 + \lp\frac{3}{80y^3}+\frac{3}{4y^2}\rp\lp\log(1+4y)-\frac{1}{\sqrt{y}}\arctan(2\sqrt{y})\rp.
\end{equation*} 
\end{thm}
\noindent
Assume that we are able to guarantee that a prime ideal $\pp$ with a fixed norm is contained in $\Ok$: this assumption provides then an explicit contribution to the estimate \eqref{RootDiscriminantLowerBoundNumerical}, which is called \textbf{local correction for the discriminant given by an ideal of norm $\Nm(\pp)$}. We denote by $C(r_1,r_2,\Nm(\pp))$ the local corrections for fields with signature $(r_1,r_2)$ given by a prime of norm $\Nm(\pp)$.

Selmane \cite{selmane1999odlyzko} computed the values of local corrections for several signatures and prime ideals: in the following tables we report the lower bounds for $|d_K|$ obtained with local corrections for fields of degree 8 and 9, in every signature, and for prime ideals $\pp$ of norm $\Nm(\pp)\leq 7$.

\begin{table}[H]
\caption{Local corrections $C(r_1,r_2,\Nm(\pp))$ for number fields of degree 8}
\begin{center}
    \begin{tabular}{l|c|c|c|c|c|r}
        $(r_1,r_2)$&(0,4) & $(2,3)$ & $(4,2)$  & $(6,1)$& (8,0) \\
        \hline
        $\Nm(\pp)=2$ & 3379343 & 11725962  & 42765027 &  163060410 & 646844001\\
        \hline
        $\Nm(\pp)=3$ & 2403757 & 8336752  & 30393063 &  115852707 & 459467465\\
        \hline
        $\Nm(\pp)=4$ & 1930702 & 6688609  & 24363884 &  92810084 & 367892401\\
        \hline
        $\Nm(\pp)=5$ & 1656110 & 5726300  & 20829049 &  79259702 & 313918560\\
        \hline
        $\Nm(\pp)=7$ & 1362891 & 4682934  & 16957023 &  64309249 & 254052210\\
    \end{tabular}
  \end{center}
  \label{localCorrection8}
  \end{table}
   \begin{table}[H]
   \caption{Local corrections $C(r_1,r_2,\Nm(\pp))$ for number fields of degree 9}
    \begin{center}
    \begin{tabular}{l|c|c|c|c|c|r}
        $(r_1,r_2)$&(1,4) & $(3,3)$ & $(5,2)$  & $(7,1)$& (9,0) \\
        \hline
        $\Nm(\pp)=2$ & 81295493 & 301476699  & 1165734091 &  4679379812 & 19422150186\\
        \hline
        $\Nm(\pp)=3$ & 57789556 & 214235371  & 828172359 &  3323651196 & 13792634200\\
        \hline
        $\Nm(\pp)=4$ & 46348899 & 171694276  & 663330644 &  2660853331 & 11037921283\\
        \hline
        $\Nm(\pp)=5$ & 39657561 & 146723910  & 566314434 &  2269968332 & 9410709985\\
        \hline
        $\Nm(\pp)=7$ & 32371189 & 119294181  & 459066389 &  1835807996 & 7596751280\\
    \end{tabular}
\end{center}
\label{localCorrection9}
\end{table}
Local corrections provide the following arithmetic consequences: if $K$ has signature $(r_1,r_2)$ and $|d_K|< C(r_1,r_2,\Nm(\pp))$, then $\Ok$ does not admit any prime ideal with norm less or equal than $\Nm(\pp)$.

This fact reflects then on the defining polynomials of the field: assume that $|d_K|\leq C(r_1,r_2,\Nm(\pp))$. If $p(x)$ is a defining polynomial of $K$ and $\alpha\in K$ is a root of $p(x)$, then we know that $|p(n)|= \Nm((\alpha-n)\Ok)$ for every $n\in\Z$, and so $p(n)$ must not be an exact multiple of every $m\in\{2,\ldots,\Nm(\pp)\}$, i.e. $m$ divides $n$ and $n/m$ is not divided by $m$.


\section{The Algorithmic Procedure}
We want to detect all the number fields $K$ with degree $n$, signature $(r_1,r_2)$ and $|d_K|\leq C(r_1,r_2,5)$, where $C(r_1,r_2,5)$ is the local correction for the signature $(r_1,r_2)$ given by a prime ideal of norm 5. We accomplish thus by constructing all the polynomials of degree $n$ having integer coefficients bounded by the values $U_m$ obtained from Theorems \ref{HunterPohstMartinet} and \ref{ThmPohst} setting $C(r_1,r_2,5)$ as upper bound of $|d_K|$. Thannks to this construction, it is clear that we are dealing with defining polynomials of HPM-elements.

The polynomials are generated ranging the values for the Newton sums $S_m$ in the intervals $[-U_m,U_m]$; from these values we create the coefficients of the polynomials with the help of the recursive relations (\ref{recursiveRelations}) and of further conditions derived from the arithmetic nature of the problem, like the fact that any evaluation of the polynomial cannot be an exact multiple of 2, 3, 4 or 5.

\begin{rmk}
As stated above, the procedure assumes that we are looking for defining polynomials of HPM-elements. There is a problem, however: unless the number field $K$ is primitive, i.e. without subfields which are not $\Q$ and $K$, nothing assures us that the defining polynomial of an HPM element $\alpha\in K$ has degree exactly equal to $n$. In fact, $\alpha$ could be contained in a proper subfield of $K$.

So we can just say that this procedure gives a complete classification only for primitive fields, which for composite degrees is still a proper subset of the considered family (though being actually a very large subset).

Fortunately, a relative version of Hunter-Pohst-Martinet's Theorem \cite{martinet1985methodes} allowed to get a complete classification of non-primitive fields up to larger upper bounds for $|d_K|$, and specifically in the following cases:
\begin{itemize}
    \item \cite{cohen1999tables} and \cite{selmane1999non} give a complete classification of non-primitive fields of degree 8 with signature $(2,3)$, $(4,2)$ and $(6,1)$ and $|d_K|\leq 6688609, 24363884$ and $92810082$ respectively;
    \item \cite{imprimitiveDegree9} gives a classification of non-primitive number fields of degree 9 with $|d_K|\leq 5 \cdot 10^7$, $4\cdot 10^9, 5 \cdot 10^9, 7 \cdot 10^9,$ $6,3 \cdot 10^{10}$ for the signatures $(1,4),(3,3), (5,2), (7,1), (9,0)$ respectively.
\end{itemize}
Thus in our procedure we can restrict ourselves just to primitive fields.
\end{rmk}
\noindent
For what concerns the algorithmic procedure, this is in fact the same we used in order to classify number fields with signature (2,3), so that we wont give all the details here, but we will just refer to what is presented in \cite{battistoniMinimum}, Section 4. In fact, we obtained Theorem \ref{TeoremaClassificazioneDiscriminanti} by following the instructions of the previous algorithm from Step 0 to Step 4.

There are nonetheless some differences which must be remarked: first of all, one should replace the previous upper bound 5762300 with the local correction $C(r_1,r_2,5)$; moreover, every feature in the algorithm related to the previously chosen degree $n=8$, like the amount of nested loops or the checks done in Step 3,  can be easily generalized for an arbitrary degree $n$.

Next, there are some additional tests that can be made already in Step 1: the polynomial $p(x)$ is kept if and only if it is constructed by Newton sums satisfying the followings restraints:  
    \begin{align}\label{PohstConditions}
&\text{If } a_1=0, \text{ then } S_3\geq 0, \nonumber\\
& S_2\geq -U_2 + \frac{2}{n}a_1^2,\nonumber\\
&|S_3|\leq \lp\frac{S_2+U_2}{2}(S_4+2(U_2-S_2)^2)\rp^{1/2}, \\
&S_4\geq -2(U_2-S_2)^2.\nonumber
\end{align} 
\noindent
The first two inequalities are proved in Cohen's book \cite{cohenAdvancedComputational}, Chapter 9. Inequality \eqref{PohstConditions} is proved by means of Cauchy-Schwartz inequality. The fourth inequality is a trivial necessary condition for the validity of the third one


.

Finally, Step 5 of the previous version is now put into Step 3, so that a candidate polynomial $p(x)$ for defining a desired number field should satisfy, together with the conditions described in Step 3, the following properties.
\begin{itemize}
    \item $p(x)$ must be an irreducible polynomial.

\item The field generated by $p(x)$ must not have prime ideals of norm less or equal than 5. This can be verified in an algorithmic way (as we explain in the next section). Moreover, the signature of $p(x)$ must be equal to $(r_1,r_2)$.

\item Given an integer $m$, define $\text{coredisc}(m)$ as the discriminant of the number field $\Q(\sqrt{m})$. Then we require $|\text{coredisc}(\text{disc}(p(x)))| < C(r_1,r_2,5).$\\
\end{itemize}
Once we have followed the instructions from Step 0 to Step 4, comprehensive of the above modifications, one just needs to replace the previous Step 5 with the following Step 5'.
\\\\
\textbf{Step 5':}  We repeat the previous steps for every value of $a_1$ between 0 and $n/2$ and for every value of $a_n$ which satisfies $|a_n|\leq (U_2/n)^{n/2}$ and is not an exact multiple of 2, 3, 4 or 5. We are left with a list of polynomials among which we select the ones generating a number field $K$ with signature $(r_1,r_2)$ with $|d_K|\leq C(r_1,r_2,5)$. 

The fields gathered are finally classified up to isomorphism and put in increasing order with respect to their absolute  discriminant.

\section{Remarks on the implementation and the results}
The theoretical ideas on which our procedure is based and the several steps composing the algorithm are very similar to what has been introduced in \cite{battistoniMinimum}, with only few differences in some of the conditions put during the tests (like the check on the size of the coredisc). The main novelty which allowed us to obtain complete tables for further signatures is the different implementation, written only in PARI/GP, which gave the following consequences and facts:
\begin{itemize}
    \item As previously mentioned, the polynomials created during the process are tested also verifying that the ideals generated by the corresponding number fields do not have norm less or equal than 5. The implementation of this process has been achieved thanks to the \textbf{ZpX-primedec()} function, written by Karim Belabas on purpose: the function is theoretically based upon the work by Ford, Pauli and Roblot (\cite{ForRob:round4}, Section 6) which use the Round 4 Algorithm in order to recover the factorization of a prime ideal from the $p$-adic factorization of a minimal polynomial of the field.

    For what concerns its efficiency, this function is an order of magnitude faster than the partial factorization given by \textbf{nfinit()} and faster than the usual decomposition function \textbf{idealprimedec()}: moreover, it is even faster whenever one deals with polynomials which give elements with small valuations for their indexes in $\Ok$. 
    \item The final check on the polynomials, suggested by Bill Allombert, concerns the size of $\text{core}(\text{disc}(p(x)))$: this test was added only some month after the signatures in degree 8 were solved. However, it allows to exclude many polynomials, because several candidate polynomials $p(x)$ have in fact core discriminants with very big size, which would force the number field discriminant to be way over the desired upper bound.\\
    The number of polynomials surviving this last condition is very small, being at most of order $10^2$, and for these one can directly compute the number field discriminant.
    \item The times of computation vary considerably and range from few hours (for signatures (2,3) and (4,2)), few days (signatures (6,1) and (1,4)), up to some months (signature (3,3)).
    \item The tables presenting all the detected number fields can be found as PARI/GP files at the website \cite{battistoniWebsite}, together with the programs written by the author, the collection of polynomials found as result of the iterations and the overview on computation times.
\end{itemize}
Finally, we present some remarks concerning the results described in Theorem \ref{TeoremaClassificazioneDiscriminanti}.
\begin{itemize}
    \item Every field in our lists is uniquely characterized by its signature and the value of its discriminant, with exceptions given only by two fields with signature $(3,3)$ and same discriminant equal to $-142989047$. These fields are given by the defining polynomials $x^9 - 4x^7 - 4x^6 + 2x^5 + 5x^4 + 6x^3 + 8x^2 + 4x + 1$ and $x^9 - 6x^7 - 9x^6 - 2x^5 + 21x^4 + 35x^3 + 23x^2 + 7x + 1$.
    \item Every field of degree 8 and every field with signature $(1,4)$ contained in our lists was already known: in fact, they are all gathered into the Kl\"{u}ners-Malle database of number fields \cite{klunersMalle}, although many of them miss from the LMFDB database \cite{lmfdb}. Our work allows to say that these are the only number fields with the corresponding signatures with discriminant less than the chosen upper bound.
    \item Concerning the fields of degree 9 and signature $(3,3)$, our procedure showed that there exist 116 such fields with $|d_K|\leq 146723910$, while the Kl\"{u}ners-Malle database only contains 62 fields of this kind. Considering the additional 54 fields, we see that 52 of them have discriminant which match with Denis Simon's table of small polynomial discriminants \cite{simonTables}. The two remaining fields satisfy instead the following properties: one of them is the field of discriminant $-142989047$ which is not isomorphic to the one given by the polynomial in Simon's list; the other one has discriminant equal to $-129079703$, which is a value not contained in Simon's lists for polynomials of degree 9 with 3 real roots, thus providing a number field and a discriminant value which were not previously known.  
    \item Every field in the list has trivial class group, and every primitive field has Galois group $G$ of the Galois closure equal to either $S_8$ or $S_9$, depending on the degree. More in detail, we have:
    \begin{itemize}
        \item[ ] 27 primitive fields out of 41 fields with signature (4,2)  (65.8\% ca.).
        \item[ ] 4 primitive fields out of 8 fields with signature (6,1)  (50\%).
        \item[ ] 63 primitive fields out of 67 fields with signature (1,4)  (94\% ca.).
        \item[ ] 112 primitive fields out of 116 fields in signature (3,3)  (96.5\% ca.).
    \end{itemize}
    Furthermore, Table \ref{TableMinimumDisc} provides the minimum values of $|d_K|$ for a field $K$ with signature $(r_1,r_2)$ and with Galois group $G$ as reported.
    \begin{table}
        \centering
        \begin{tabular}{|c|c|c|c|}
            $n$ & $(r_1,r_2)$ & $G$ & minimum $|d_K|$\\
            \hline
            8 & (4,2) & $S_8$ &15908237\\
            \hline
            8 & (6,1) & $S_8$ &65106259\\
            \hline
            9 & (1,4) & $S_9$ &29510281\\
            \hline
            9 & (3,3) & $S_9$ &109880167\\
            \hline
        \end{tabular}
        \caption{Minimum discriminants for specific Galois groups}
        \label{TableMinimumDisc}
    \end{table}
    \item Although the algorithm classifies only primitive fields, every non-primitive field with $|d_K|\leq C(n,r_1,5)$ was displayed as output. 
    \item The groups in \cite{battistoniWebsite} are presented according to the LMFDB notation: every group is denoted by $n$T$q$, where $n$ is the degree of the corresponding field and $q$ is the label of the group as transitive subgroup of $S_n$: the choice of the label is based upon Hulpke's algorithm for the classification of transitive subgroups of $S_n$ \cite{hulpkeTransitive}. If the group has an easy form, like the dihedral group $D_n$ or the symmetric group $S_n$, then the classical name of the group is written together with the LMFDB label.

\end{itemize}


\begin{thebibliography}{10}

\bibitem{regulators}
S.~Astudillo, F.~Diaz~y Diaz, and E.~Friedman.
\newblock Sharp lower bounds for regulators of small-degree number fields.
\newblock {\em J. Number Theory}, 167:232--258, 2016.

\bibitem{battistoniWebsite}
F.~Battistoni.
\newblock {\em {Tables of Number Fields}}.
\newblock available at
  \url{http://www.mat.unimi.it/users/battistoni/index.html}.

\bibitem{battistoniMinimum}
F.~Battistoni.
\newblock The minimum discriminant of number fields of degree 8 and signature
  {$(2,3)$}.
\newblock {\em J. Number Theory}, 198:386--395, 2019.

\bibitem{belabasCubic}
K.~Belabas.
\newblock A fast algorithm to compute cubic fields.
\newblock {\em Math. Comp.}, 66(219):1213--1237, 1997.

\bibitem{bergeOlivierMartinet}
A.-M. Berg\'{e}, J.~Martinet, and M.~Olivier.
\newblock The computation of sextic fields with a quadratic subfield.
\newblock {\em Math. Comp.}, 54(190):869--884, 1990.

\bibitem{QuarticFields}
J.~Buchmann, D.~Ford, and M.~Pohst.
\newblock Enumeration of quartic fields of small discriminant.
\newblock {\em Math. Comp.}, 61(204):873--879, 1993.

\bibitem{cohenAdvancedComputational}
H.~Cohen.
\newblock {\em Advanced topics in computational number theory}, volume 193 of
  {\em Graduate Texts in Mathematics}.
\newblock Springer-Verlag, New York, 2000.

\bibitem{cohen1999tables}
H.~Cohen, F.~Diaz~y Diaz, and M.~Olivier.
\newblock Tables of octic fields with a quartic subfield.
\newblock {\em Math. Comp.}, 68(228):1701--1716, 1999.

\bibitem{davenportHeilbronn}
H.~Davenport and H.~Heilbronn.
\newblock On the density of discriminants of cubic fields. {II}.
\newblock {\em Proc. Roy. Soc. London Ser. A}, 322(1551):405--420, 1971.

\bibitem{y1980tables}
F.~Diaz~y Diaz.
\newblock {\em Tables minorant la racine {$n$}-i\`eme du discriminant d'un
  corps de degr\'e {$n$}}, volume~6 of {\em Publications Math\'ematiques
  d'Orsay 80 [Mathematical Publications of Orsay 80]}.
\newblock Universit\'e de Paris-Sud, D\'epartement de Math\'ematique, Orsay,
  1980.

\bibitem{diazSepticUneSeulePlace}
F.~Diaz~y Diaz.
\newblock Valeurs minima du discriminant des corps de degr\'{e} {$7$} ayant une
  seule place r\'{e}elle.
\newblock {\em C. R. Acad. Sci. Paris S\'{e}r. I Math.}, 296(3):137--139, 1983.

\bibitem{diazSepticCertain}
F.~Diaz~y Diaz.
\newblock Valeurs minima du discriminant pour certains types de corps de
  degr\'{e} {$7$}.
\newblock {\em Ann. Inst. Fourier (Grenoble)}, 34(3):29--38, 1984.

\bibitem{y1987petits}
F.~Diaz~y Diaz.
\newblock Petits discriminants des corps de nombres totalement imaginaires de
  degr\'e {$8$}.
\newblock {\em J. Number Theory}, 25(1):34--52, 1987.

\bibitem{diazSeptic}
F.~Diaz~y Diaz.
\newblock Discriminant minimal et petits discriminants des corps de nombres de
  degr\'{e} {$7$} avec cinq places r\'{e}elles.
\newblock {\em J. London Math. Soc. (2)}, 38(1):33--46, 1988.

\bibitem{imprimitiveDegree9}
F.~Diaz~y Diaz and M.~Olivier.
\newblock Imprimitive ninth-degree number fields with small discriminants.
\newblock {\em Math. Comp.}, 64(209):305--321, 1995.
\newblock With microfiche supplement.

\bibitem{ForRob:round4}
D.~Ford, S.~Pauli, and X.F. Roblot.
\newblock A fast algorithm for polynomial factorization over {$\mathbb{Q}\sb
  p$}.
\newblock {\em J. Th\'eor. Nombres Bordeaux}, 14(1):151--169, 2002.

\bibitem{hulpkeTransitive}
A.~Hulpke.
\newblock Constructing transitive permutation groups.
\newblock {\em J. Symbolic Comput.}, 39(1):1--30, 2005.

\bibitem{klunersMalle}
J.~Kl\"{u}ners and G.~Malle.
\newblock A database for number fields.
\newblock Available at \url{http://galoisdb.math.upb.de/home}.

\bibitem{martinet1985methodes}
J.~Martinet.
\newblock Methodes g\'eom\'etriques dans la recherche des petits discriminants.
\newblock In {\em S\'eminaire de th\'eorie des nombres, {P}aris 1983--84},
  volume~59 of {\em Progr. Math.}, pages 147--179, 1985.

\bibitem{odlyzko1990bounds}
A.~M. Odlyzko.
\newblock Bounds for discriminants and related estimates for class numbers,
  regulators and zeros of zeta functions: a survey of recent results.
\newblock {\em S\'em. Th\'eor. Nombres Bordeaux (2)}, 2(1):119--141, 1990.

\bibitem{olivierSexticQuadratic1}
M.~Olivier.
\newblock Corps sextiques contenant un corps quadratique. {I}.
\newblock {\em S\'{e}m. Th\'{e}or. Nombres Bordeaux (2)}, 1(1):205--250, 1989.

\bibitem{olivierQuadratic2}
M.~Olivier.
\newblock Corps sextiques contenant un corps quadratique. {II}.
\newblock {\em S\'{e}m. Th\'{e}or. Nombres Bordeaux (2)}, 2(1):49--102, 1990.

\bibitem{olivierPrimitifs}
M.~Olivier.
\newblock Corps sextiques primitifs.
\newblock {\em Ann. Inst. Fourier (Grenoble)}, 40(4):757--767 (1991), 1990.

\bibitem{olivierSextiquesCubique}
M.~Olivier.
\newblock Corps sextiques contenant un corps cubique. {III}.
\newblock {\em S\'{e}m. Th\'{e}or. Nombres Bordeaux (2)}, 3(1):201--245, 1991.

\bibitem{olivierSextiquesPrimitifs}
M.~Olivier.
\newblock Corps sextiques primitifs. {IV}.
\newblock {\em S\'{e}m. Th\'{e}or. Nombres Bordeaux (2)}, 3(2):381--404, 1991.

\bibitem{olivierSexticCubic1992}
M.~Olivier.
\newblock The computation of sextic fields with a cubic subfield and no
  quadratic subfield.
\newblock {\em Math. Comp.}, 58(197):419--432, 1992.

\bibitem{pohstDegree7TotallyReal}
M.~Pohst.
\newblock The minimum discriminant of seventh degree totally real algebraic
  number fields.
\newblock In {\em Number theory and algebra}, pages 235--240. 1977.

\bibitem{pohst1982computation}
M.~Pohst.
\newblock On the computation of number fields of small discriminants including
  the minimum discriminants of sixth degree fields.
\newblock {\em J. Number Theory}, 14(1):99--117, 1982.

\bibitem{pohst1990minimum}
M.~Pohst, J.~Martinet, and F.~Diaz~y Diaz.
\newblock The minimum discriminant of totally real octic fields.
\newblock {\em J. Number Theory}, 36(2):145--159, 1990.

\bibitem{poitou1977petits}
G.~Poitou.
\newblock Sur les petits discriminants.
\newblock In {\em S\'eminaire {D}elange-{P}isot-{P}oitou, 18e ann\'ee:
  (1976/77), {T}h\'eorie des nombres, {F}asc. 1 ({F}rench)}, pages Exp. No. 6,
  18. Secr\'etariat Math., Paris, 1977.

\bibitem{QuinticFields}
A.~Schwarz, M.~Pohst, and F.~Diaz~y Diaz.
\newblock A table of quintic number fields.
\newblock {\em Math. Comp.}, 63(207):361--376, 1994.

\bibitem{selmane1999non}
S.~Selmane.
\newblock Non-primitive number fields of degree eight and of signature
  {$(2,3)$}, {$(4,2)$} and {$(6,1)$} with small discriminant.
\newblock {\em Math. Comp.}, 68(225):333--344, 1999.

\bibitem{selmane1999odlyzko}
S.~Selmane.
\newblock Odlyzko-{P}oitou-{S}erre lower bounds for discriminants for some
  number fields.
\newblock {\em Maghreb Math. Rev.}, 8(1-2):151--162, 1999.

\bibitem{serre1986minorations}
{J}.~{P}. {S}erre.
\newblock Minorations de discriminants, note of october 1975, published on pp.
  240-243 in vol. 3 of {J}ean-{P}ierre {S}erre, collected papers, 1986.

\bibitem{simonTables}
D.~Simon.
\newblock Petits discriminants de polynomes irréductibles.
\newblock available at
  \url{https://simond.users.lmno.cnrs.fr/maths/TableSmallDisc.html}.

\bibitem{takeuchiNine}
K.~Takeuchi.
\newblock Totally real algebraic number fields of degree 9 with small
  discriminant.
\newblock {\em Saitama Math. J.}, 17:63--85 (2000), 1999.

\bibitem{lmfdb}
{The LMFDB Collaboration}.
\newblock The l-functions and modular forms database.
\newblock available at \url{http://www.lmfdb.org}, 2013.

\bibitem{pari}
{The PARI~Group}, Univ. Bordeaux.
\newblock {\em {PARI/GP version {\tt 2.11.0}}}, 2018.
\newblock available at \url{http://pari.math.u-bordeaux.fr/}.

\end{thebibliography}
\end{document}